\documentclass{article}
\usepackage{amsmath}
\usepackage{amssymb}
\usepackage[dvips]{graphics}

\def\v{\varepsilon}

\def\t{\theta}

\def\a{\alpha}
\def\b{\beta}

\def\di{\displaystyle}

\begin{document}
\title{\bf Stability of viscous shock wave for compressible Navier-Stokes equations with free boundary}
\author{\bf Feimin Huang$^\dag$, Xiaoding Shi$^{\dag\dag}$, Yi Wang$^\dag$}
\date{} \maketitle
\noindent{\small  $^\dag$Institute of Applied Mathematics, AMSS,
Academia Sinica, Beijing 100190, China

\noindent $^{\dag\dag}$Department of Mathematics, Graduate School
of Science, Beijing University of Technology and Chemical, Beijing
100029, China}

\

 \noindent {\bf Abstract:} A free boundary problem for the one-dimensional compressible Navier-Stokes equations is investigated.
 The asymptotic stability of the viscous
 shock wave is established under some smallness conditions. The proof is given
 by an elementary energy estimate.

\section{Introduction}
We consider the system of  viscous and heat conductive fluid in
the Eulerian coordinate
$$
\left\{
\begin{array}{llll}
{\displaystyle \rho_{t}+(\rho u)_{\tilde{x}}=0,}
 \\
{\displaystyle (\rho u)_{t}+(\rho u^{2}+p)_{\tilde{x}}=\mu
u_{\tilde{x}\tilde{x}},} \\
{\displaystyle \left[\rho(e+\frac{u^{2}}{2})\right]_{t}+
\left[\rho u(e+\frac{u^{2}}{2})+pu\right]_{\tilde{x}}
=\kappa\theta_{\tilde{x}\tilde{x}}+ (\mu
uu_{\tilde{x}})_{\tilde{x}},}
\end{array}
\right.
\eqno(1.1)
$$
  where $u(\tilde{x},t)$ is the velocity,
$\rho(\tilde{x},t)>0$ is the density,
 $\theta(\tilde{x},t)$ is the absolute temperature,
 $p=p(\rho,\theta)$ is the pressure, $e=e(\rho,\theta)$ is
the internal energy, $\mu>0$ is the viscosity constant, $\kappa>0$
is the coefficient of heat conduction. Here we consider the
perfect gas case, that is
$$
p=R\theta\rho,\quad e=\frac{R\theta}{\gamma-1}+const., \eqno(1.2)
$$
where $\gamma>1$ is the adiabatic constant and $R>0$ the gas
constant.

There has been a large literature on the asymptotic behaviors of the
solutions to the system (1.1). However, most results are concerned
with the initial value problem. We refer to \cite{[5]}-\cite{[8]},
\cite{[10]}-\cite{[11]} and references therein. Recently the initial
boundary value problem (IBVP) attracts an increasing interest
because it has more physically meanings and of course produces new
mathematical difficulty due to the boundary effect. We refer to
\cite{[2]}, \cite{[9]}, \cite{[12]}, \cite{[13]} for $2\times2$ case
and \cite{[3]}, \cite{[4]}, \cite{[15]} for $3\times3$ case.
However, there is few result on the asymptotic stability of the
viscous shock wave to IBVP of the full compressible Navier-stokes
equation (1.1) due to various difficulties. Therefore, the
asymptotic stability of the viscous shock wave to IBVP for (1.1) is
our main purpose of the present paper. We shall consider a free
boundary problem of the full compressible Navier-Stokes equations
whose boundary conditions read
$$
\left\{
\begin{array}{llll}
{\displaystyle (p-\mu
u_{\tilde{x}})\bigl|_{\tilde{x}=\tilde{x}(t)}=p_{0},}\\
{\displaystyle \theta|_{\tilde{x}=\tilde{x}(t)}=\theta_{-}>0,}\\
{\displaystyle
\frac{d\tilde{x}(t)}{dt}=\tilde{u}(\tilde{x}(t),t),\
\tilde{x}(0)=0,\ t>0,}
\end{array}
\right.
\eqno (1.3)
$$
and initial
data
$$
(\rho,u,\theta)\bigl|_{t=0}=(\rho_{0},u_{0},\theta_{0})(x)\rightarrow(\rho_{+},u_{+},\theta_{+})\
\mathrm{as}\ \tilde{x}\rightarrow+\infty, \eqno(1.4)
$$ where
$p_{0}>0,\t_->0,\rho_{+}>0, \theta_{+}>0, u_{+}$ are prescribed
constants. Here the boundary condition (1.3) means the gas is
attached at the boundary $\tilde{x}=\tilde{x}(t)$ to the atmosphere
with pressure $p_{0}$(see \cite{[13]}). We of course assume  the
initial data satisfy the boundary condition as compatibility
condition.

 Since the
boundary condition (1.3) means the particles always stay on the
free boundary $\tilde{x}=\tilde{x}(t)$, if we use the Lagrangian
coordinates, then the free boundary becomes a fixed boundary. Thus
we transform the Eulerian coordinates $(x,t)$ by
$$x=\int_{\tilde{x}(t)}^{\tilde{x}}\rho(y,t)dy,\ t=t,$$
and then change the free boundary value problem (1.1)-(1.4) into
$$\left\{
\begin{array}{ll}
\displaystyle v_{t}-u_{x}=0,, &x>0,t>0,\\
\displaystyle u_{t}+p_{x}=\mu(\frac {u_{x}}v)_{x},
 &x>0,t>0,\\
\displaystyle\bigl(e+\frac{u^{2}}{2}\bigr)_{t}+
(pu)_{x}=(\kappa\frac{\theta_{x}}{v}+\mu
\frac{uu_{x}}{v})_{x},&x>0,t>0,\\[3mm]
\displaystyle (p-\mu \frac{u_{x}}{v})|_{x=0}=p_{0},\qquad
  \theta|_{x=0}=\theta_{-},&\\
\displaystyle(v,u,\theta)(x,0)=(v_{0},u_{0},\theta_{0})(x)\rightarrow(v_{+},u_{+},\theta_{+})&\
\mathrm{as} \ x\rightarrow+\infty ,
\end{array}
\right. \eqno(1.5)$$ where $v=\frac{1}{\rho}$ is the specific
volume. Since the domain we consider here in the Lagrange
coordinates is $\{x>0,t>0\}$, we only need to consider the
stability of the 3-viscous shock wave.

Before formulating our main result, we briefly recall some results
of the shock wave for the  inviscid system of (1.1).
 That is, we consider the system (1.5) without viscosity
$$\left\{
\begin{array}{llll}
{\displaystyle v_{t}-u_{x}=0,}\\
{\displaystyle u_{t}+p_{x}=0,}\\
{\displaystyle\bigl(e+\frac{u^{2}}{2}\bigr)_{t}+ (pu)_{x}=0,}
 \end{array}
\right. \eqno(1.6)$$ with the Riemann initial data
$$(v_0,u_0,\theta_0)(x)=\left\{
\begin{array}{llll}(v_{-},u_{-},\theta_{-}),\ x>0,\\
(v_{+},u_{+},\theta_{+}),\ x<0.\end{array}\right.\eqno(1.7)$$

It is well known (for example, see \cite{[14]}) that the Riemann
problem (1.6)-(1.7) admits a 3-shock wave if and only if the two
states $(v_{\pm},u_{\pm},\theta_{\pm})$ satisfy the so-called
 Rankine-Hugoniot condition
$$
\left\{
\begin{array}{llll}
{\displaystyle-s(v_{+}-v_{-})-(u_{+}-u_{-})=0},\\
 {\displaystyle
-s (u_{+}-u_{-})+(p_{+}-p_{-})=0,}\\
{\displaystyle-s\left[(e_{+}+\frac{u_{+}^{2}}{2})-
(e_{-}+\frac{u_{-}^{2}}{2})\right]+(p_{+}u_{+}-p_{-}u_{-})=0},
\end{array}\right.\eqno(1.8)
$$
and the Lax's entropy condition
$$0<\lambda_3^+<s<\lambda_3^-,\eqno(1.9)$$
where $p_{\pm}=p(v_{\pm},\theta_{\pm}),
e_{\pm}=e(v_{\pm},\theta_{\pm})$ and $\lambda_3=\frac{\sqrt{\gamma R
\theta}}{v}$ is the third eigenvalue of the inviscid system (1.6).
And the shock speed $s$ is uniquely determined by
$(v_{\pm},u_{\pm},\theta_{\pm})$ with (1.8). If the right state
$(v_{+},u_{+},\theta_{+})$ is given, it is easy to know that there
exists a 3-shock curve $S_{3}(v_{+},u_{+},\theta_{+})$ starting from
$(v_{+},u_{+},\theta_{+})$. For any point $(v,u,\theta)\in
S_{3}(v_{+},u_{+},\theta_{+})$, there exists a unique 3-shock wave
connecting $(v,u,\theta)$ with $(v_{+},u_{+},\theta_{+})$. Our
assumptions on the boundary values are

\

(A1).  Let $(v_{+},u_{+},\theta_{+})$ and $\theta_{-}$ be given,
there exist unique $v_{-},u_{-}$ such that
$(v_{-},u_{-},\theta_{-})\in S_{3}(v_{+},u_{+},\theta_{+})$.

\

(A2).  $p_{0}=\frac{R\theta_{-}}{v_{-}}:=p_{-}.$

\

Remark1. The assumption (A1) is natural.

Remark2.  The condition (A2) means that we only consider the
stability of a single viscous shock wave.

\

It is known that the system (1.5) admits smooth travelling wave
solution with shock profile $(V,U,\Theta)(x-s t)$ under the
conditions (1.8) and (1.9) (see \cite{[1]}). Such travelling wave
has been shown nonlinear stable for the initial value problem, see
\cite{[5]} and \cite{[7]}. A natural question is whether the
travelling wave is stable or not for the initial boundary value
problem. In this paper, we give a positive answer for the free
boundary problem (1.1)-(1.4) or (1.5). Our main result is, roughly
speaking, as follows. The precise statement is given in theorem 2.1
below.

\

 Let $(v_{+},u_{+},\theta_{+})$ and $\theta_{-}$ be given and the
 assumptions (A1) and (A2) hold, then the 3-viscous shock wave
 connecting $(v_{-},u_{-},\theta_{-})$
with $(v_{+},u_{+},\theta_{+})$ is asymptotically stable.

\

The plan of this paper is as follows. After stating the notations,
in section 2, we give some properties of the viscous shock wave
and the main Theorem 2.1. In Section 3, we reformulate the
original problem to a new initial boundary value problem. The
proof of the Theorem 2.1 is given in section 4 by the elementary
energy method. In section 5, we prove the local existence of the
solution by the iteration method.

\

\noindent\textbf{Notation:}  Throughout this paper, several
positive generic constants which are independent of $T,\beta$ and
$\alpha$ are denoted by $C$ without confusions. For function
spaces, $H^l(\mathbb{R}^{+})$ denotes the $l$-th order Sobolev
space with its norm
$$
\|f\|_l=(\sum^l_{j=0}\|\partial^j_xf\|^2)^\frac{1}{2}, \quad {\rm
when}~\|\cdot\|:=\|\cdot\|_{L^2(\mathbb{R}^{+})}.\eqno(1.10)
$$

\section{Preliminaries and Main Result}\label{sec:int}

 We first recall some properties of the 3-viscous
 shock wave. The shock profile $(V,U,\Theta)(\xi), \xi=x-st$, is
 determined by
 $$
 \left\{
\begin{array}{lll}
{\displaystyle -s V'-U'=0,}\\
{\displaystyle -s
U'+P'=\mu\left(\frac{U'}{V}\right)',}\\
{\displaystyle
-s\left(E+\frac{U^{2}}{2}\right)'+\left(PU\right)'=\left(\kappa\frac{\Theta'}{V}+\mu
\frac{UU'}{V}\right)',}\\
{\displaystyle\left(V,U,\Theta\right)(\pm\infty)=(v_{\pm},u_{\pm},\theta_{\pm}),}
\end{array}
\right.
\eqno(2.1)
$$
where $P=R\Theta/V$, $E=R\Theta/(\gamma-1)+const.$,
$(v_{\pm},u_{\pm},\theta_{\pm})$ satisfy R-H condition (1.8) and
entropy condition (1.9) and $s$ is determined by (1.8).
Integrating (2.1) on $(-\infty,\xi)$ gives
$$
 \left\{
\begin{array}{lll}
{\displaystyle \frac{s \mu
V_{\xi}}{V}=-\left[P+s^{2}(V-\frac{b_{1}}{s^{2}})\right],}\\
{\displaystyle \frac{\kappa\Theta_{\xi}}{s
V}=-\left[E-\frac{s^{2}(V-\frac{b_{1}}{s^{2}})^{2}}{2}+\frac{b_{1}^{2}}{2s^{2}}
-b_{2}\right],}\\
U=-(s V+a),
\end{array}
\right.
\eqno (2.2)
$$ where
$p_{\pm}=R\theta_{\pm}/v_{\pm}$,
$e_{\pm}=R\theta_{\pm}/(\gamma-1)+const.$, $a=-(sv_\pm+u_\pm)$,
$b_{1}=p_{\pm}+s^{2}v_{\pm}$ and
$b_{2}=e_{\pm}+p_{\pm}v_{\pm}+s^{2}v_{\pm}^{2}/2.$ From \cite{[1]}
and \cite{[7]}, we have the following proposition:

\

\noindent\textbf{Proposition 2.1.} {\it Assume that the two states
$(v_\pm,u_\pm,\theta_\pm)$ satisfy the conditions (1.8) and (1.9),
then there exists a unique shock profile $(V,U,\Theta)(\xi)$, up
to a shift, of system (2.1). Moreover, there are  positive
constants $c_1$ and $c_2$ independent of $\gamma>1$ such that for
$\xi\in\mathbb{R}$,
$$
\left\{
\begin{array}{lll}
{\displaystyle s V_{\xi}=-U_{\xi}>0,\
s\Theta_{\xi}<0,(|V-v_{\pm}|,|U-u_{\pm}|)\le c_1d e^{-c_2d|\xi|}}\\
{\displaystyle |\Theta-\theta_{\pm}|\le c_1(\gamma-1)d
e^{-c_2d|\xi|},(|V_{\xi}|,|V_{\xi\xi}|,|\Theta_{\xi\xi}|)\leq
c_1d^2 e^{-c_2d|\xi|},}\\
{\displaystyle|\Theta_{\xi}|\leq c_1(\gamma-1)d^2e^{-c_2d|\xi|},\
|\frac{\Theta_{\xi}}{V_{\xi}}|\leq c_1(\gamma-1),}\\
{\displaystyle s^2=\frac{\gamma R\theta_-(1-d_1)}{v_+v_-},
d_1=\frac{d_2}{1+d_2},
d_2=\frac{(\gamma-1)d}{2v_+},}
\end{array}
\right.
\eqno(2.3)
$$ where
$d=v_+-v_-.$}

\

As pointed out by Liu \cite{[5]}, a generic perturbation of viscous
shock wave produces not only a shift $\alpha$ but also diffusion
waves, which decay to zero with a rate $(1+t)^{-\frac12}$, for the
Cauchy problem. That is the solution of the compressible
Navier-Stokes equations asymptotically tends to the translated
travelling wave $(V,U,\Theta)(x-st+\alpha)$. The shift $\alpha$ is
explicitly determined by the initial value. Similar to the Cauchy
problem, the shift $\alpha$ is also expected for IBVP. For a kind of
initial boundary value problem, in which the velocity is zero on the
boundary, Matsumura and Mei \cite{[9]} developed a new way to
determine the shift $\alpha$. A byproduct of \cite{[9]} showed that,
unlike the Cauchy problem, there is no diffusion wave for IBVP due
to the boundary effect. This new idea has been used by many authors
to treat the initial boundary value problem of the system (1.5) or
other related systems (see \cite{[2]}, \cite{[12]}, \cite{[13]}). In
the spirit of \cite{[9]}, we calculate the shift $\alpha$ for the
IBVP (1.5).

We consider the situation where the initial data
$(v_{0},u_{0},\theta_{0})$ are given in a neighborhood of
$(V,U,\Theta)(x-\beta)$ for some large constant $\beta>0$. That
is, we require the viscous shock wave is far from the boundary
initially. Here we can not directly apply the idea of \cite{[9]}
to compute the shift $\alpha$ since the velocity $u(0,t)$ on the
boundary is not given,while in \cite{[9]}, the velocity is zero on
the boundary and the conservation of the mass $(1.5)_1$ is then
used to determine the shift $\alpha$. instead of $(1.5)_1$, we use
the conservation of momentum $(1.5)_2$ to determine the shift
$\alpha$
 because $p-\mu\frac{u_x}{v}$ is given on the boundary for the
 IBVP (1.5). From $(1.5)_2$ and $(2.1)_2$, we have
$$
(u-U)_{t}=-[p(v,\theta)-P(V,\Theta)]_{x}+\mu\left(\frac{u_{x}}{v}\right)_{x}-\mu\left(\frac{U_{x}}{V}\right)_{x},
\eqno(2.4)
$$
where $(V,U)=(V,U)(x-s t+\alpha-\beta)$. Integrating
(2.4) over $[0,\infty)$ with respect to $x$ and using (2.1) and
(A2) yield
$$
\begin{array}{lll}
{\displaystyle\frac{d}{dt}\int_{0}^{+\infty}[u(x,t)-U(x-s
t+\alpha-\beta)]dx}\\
{\displaystyle=p_{-}-P(V,\Theta)(-s
t+\alpha-\beta)+\mu\frac{U'}{V}(-s
t+\alpha-\beta)}\\
{\displaystyle=-s (U(-s t+\alpha-\beta)-u_{-}).}
\end{array}
\eqno(2.5)
$$
Integrating (2.5) with respect to
$t$, we have
$$
\begin{array}{lll}
{\displaystyle\int_{0}^{+\infty}[u(x,t)-U(x-s
t+\alpha-\beta)]dx}\\
{\displaystyle=\int_{0}^{+\infty}[u_{0}-U(x+\alpha-\beta)]dx-\int_{0}^{t}s
(U(-s \tau+\alpha-\beta)-u_{-})d\tau.}
\end{array}
\eqno(2.6)
$$ We
define
$$
\begin{array}{ll}
I(\alpha):=&{\displaystyle
\int_{0}^{+\infty}[u_{0}-U(x+\alpha-\beta)]dx}\\
&{\displaystyle-\int_{0}^{+\infty}s \left(U(-s
t+\alpha-\beta)-u_{-}\right)dt.}
\end{array}
\eqno(2.7)
$$
It follows that
$$
\begin{array}{ll}
I'(\alpha)=&{\displaystyle
-\int_{0}^{+\infty}U'(x+\alpha-\beta)dx-s\int_{0}^{\infty}U'(-s
\tau+\alpha-\beta)]d\tau}\\
 &{\displaystyle = u_{-}-u_{+}.}\end{array}\eqno(2.8)$$
Expectation $\lim_{t\to \infty}\int_{0}^{+\infty}[u(x,t)-U(x-s
t+\alpha-\beta)]dx=I(\alpha)=0$ gives
$$\alpha=\frac{1}{u_{+}-u_{-}}I(0),\eqno (2.9)$$ and
$$\begin{array}{lll}
{\displaystyle \int_{0}^{+\infty}[u(x,t)-U(x-s
t+\alpha-\beta)]dx}\\
{\displaystyle= s\int_{t}^{+\infty}[U(-s
\tau+\alpha-\beta)-u_-]d\tau\leq c_1e^{-c_2d|-st+\alpha-\beta|}\
\mathrm{as}\ t\rightarrow+\infty.}
\end{array}
\eqno(2.10)
$$
Therefore the shift$\alpha$ is uniquely determined  by the initial
value.
\

To state our main theorem, we suppose that for some $\beta>0$
$$
\left(v_{0}(x)-V(x-\beta),u_{0}(x)-U(x-\beta),\theta_{0}(x)-\Theta(x-\beta)\right)\in
H^{1}\cap L^{1}.\eqno (2.11)
$$
Let
$$
\begin{array}{lll}
{\displaystyle
(\widetilde{\Phi}_{0},\widetilde{\Psi}_{0})(x)=-\int_{x}^{+\infty}\left[v_{0}(y)-V(y-\beta),
u_{0}(y)-U(y-\beta)\right]dy,}\\
{\displaystyle \widetilde{W}_{0}(x)
=-\int_{x}^{+\infty}\left[(e_{0}+\frac{u_{0}^{2}}{2})(y)-(E+\frac{U^{2}}{2})(y-\beta)\right]dy.}
\end{array}
\eqno(2.12)
$$
Assume that
$$
(\widetilde{\Phi}_{0},\widetilde{\Psi}_{0},\widetilde{W}_{0})\in L^{2}.
\eqno (2.13)
$$
Our main result is

\

\noindent \textbf{Theorem 2.1.}  {\it Suppose that the assumptions
(A1) and (A2) hold. Let
 $(V,U,\Theta)(\xi)$ be the travelling wave solution satisfying (2.1).
 Assume that $1<\gamma\leq 2$ and
 (2.11-2.13) hold,
then there exists positive
 constants $\delta_0$ and $\varepsilon_0$
 such that if  $$(\gamma-1)d\leq\delta_0,\eqno (2.14)$$
and
 $$
 \|(\widetilde{\Phi}_{0},\widetilde{\Psi}_{0},\frac{\widetilde{W}_{0}}{\sqrt{\gamma-1}})\|_2
 +e^{-c_2d\beta}\leq\varepsilon_0,
 \eqno(2.15)
 $$
then the system (1.5) has a unique global solution
 $(v,u,\theta)(x,t)$ satisfying
 $$
 \begin{array}{lll}
v(x,t)-V(x-s t+\alpha-\beta)\in
 C([0,\infty),H^{1})\cap L^{2}(0,\infty;H^{1}),\\
 u(x,t)-U(x-s t+\alpha-\beta)\in
 C([0,\infty),H^{1})\cap L^{2}(0,\infty;H^{2}),\\
 \theta(x,t)-\Theta(x-s t+\alpha-\beta)\in
 C([0,\infty),H^{1})\cap L^{2}(0,\infty;H^{2}),\\
 \end{array}
 \eqno(2.16)
 $$
and
 $$
 \sup_{x\in\mathbb{R}_{+}}\bigl|(v,u,\theta)(x,t)-(V,U,\Theta)(x-st+\alpha-\beta)
 \bigr|\longrightarrow 0,\mathrm{\ as\ }t\rightarrow+\infty,
 \eqno(2.17)
 $$
 where $\alpha=\alpha(\beta)$ is determined by (2.9).}

\section{Reformulation of the Original Problem}\label{sec:int}

\noindent Let
$$
(v,u,\theta)(x,t)=(V,U,\Theta)(x-s t+\alpha-\beta)+(\phi,\psi,w)(x,t),
\eqno (3.1)
$$
then we rewrite the system (1.5) as
$$
\left\{
\begin{array}{llll}
{\displaystyle \phi_{t}-\psi_{x}=0,}\\
{\displaystyle
\psi_{t}+R\left(\frac{\Theta+w}{V+\phi}-\frac{\Theta}{V}\right)_{x}
=\mu\left[\frac{\psi_{x}}{V+\phi}+\left(\frac{1}{V+\phi}-
\frac{1}{V}\right)U_{x}\right]_{x},}\\[3mm]
{\displaystyle\left(\frac{R}{\gamma-1}w+\frac{\psi^{2}}{2}+U\psi\right)_{t}
+R\left[\frac{\Theta+w}{V+\phi}\psi+(\frac{\Theta+w}{V+\phi}-\frac{\Theta}{V})U\right]_{x}}\\[3mm]
{\displaystyle\qquad
=\kappa\left[\frac{w_{x}}{V+\phi}+(\frac{1}{V+\phi}-\frac{1}{V})\Theta_{x}\right]_{x}
}\\[3mm]
{\displaystyle\qquad+\mu\left[\frac{\psi\psi_{x}+U\psi_{x}+U_{x}\psi}{V+\phi}+
(\frac{1}{V+\phi}-\frac{1}{V})
UU_{x}\right]_{x},}\\[3mm]
{\displaystyle
w|_{x=0}=\theta_--\Theta(-st+\alpha-\beta),\quad\left(\frac{R\theta_{-}}{V+\phi}-
\mu\frac{U_{x}+\psi_{x}}{V+\phi}\right)\bigr|_{x=0}=p_{-},}\\[3mm]
 {\displaystyle(\phi,\psi,w)|_{t=0}=(\phi,\psi,w)(x,0):=
(\phi_{0},\psi_{0},w_{0})(x).}
\end{array}
\right.
\eqno(3.2)
$$
We define
$$
\begin{array}{lll}
{\displaystyle(\Phi,\Psi)(x,t)=-\int_{x}^{+\infty}\left(\phi,\psi
\right)(y,t)dy,}\\
{\displaystyle
W(x,t)=-\int_{x}^{+\infty}\left(e+\frac{u^{2}}{2}\right)(y,t)
-\left(E+\frac{U^{2}}{2}\right)(y-s t
+\alpha-\beta)dy.}
\end{array}
\eqno(3.3)
$$
Then we have
$$
(\phi,\psi,w)=\left(\Phi_{x},\Psi_{x},\frac{\gamma-1}{R}[W_{x}-
(\frac{1}{2}\Psi_{x}^{2}+U\Psi_{x})]\right).
\eqno (3.4)
$$
Integrating (3.2) with respect to $x$ yields
$$
\left\{
\begin{array}{lll}
{\displaystyle \Phi_{t}-\Psi_{x}=0,}\\
{\displaystyle\Psi_{t}+R\left(\frac{\Theta+w}{V+\Phi_{x}}-\frac{\Theta}{V}\right)
=\frac{\mu\Psi_{xx}}{V+\Phi_{x}}+\left(\frac{\mu}{V+\Phi_{x}}-\frac{\mu}{V}\right)U_{x},}\\
{\displaystyle
W_{t}+R\left(\frac{\Theta+w}{V+\Phi_{x}}-\frac{\Theta}{V}\right)U
+R\frac{\Theta+w}{V+\Phi_{x}}\Psi_{x}}\\
{\displaystyle\quad=\frac{\kappa w_{x}}{V+\Phi_{x}}
+\left(\frac{\kappa}{V+\Phi_{x}}-\frac{\kappa}{V}\right)\Theta_{x}}\\
{\displaystyle\quad+
\frac{\mu}{V+\Phi_{x}}(\Psi_{x}\Psi_{xx}+U_{x}\Psi_{x}+U\Psi_{xx})+
(\frac{\mu}{V+\Phi_{x}}-\frac{\mu}{V}) UU_{x}}.
\end{array}
\right.
\eqno(3.5)
$$
Introduce a new variable
$$
\widehat{W}=\frac{\gamma-1}{R}(W-U\Psi),
\eqno(3.6)
$$
then we write $w$ in the form
$$
w=\widehat{W}_{x}+\frac{\gamma-1}{R}\left(U_{x}\Psi-\frac{\Psi_{x}^{2}}{2}\right),
\eqno(3.7)
$$
and transform the system (3.5) into
$$
\left\{
\begin{array}{lll}
{\displaystyle \Phi_{t}-\Psi_{x}=0,}\\
{\displaystyle\Psi_{t}-\frac{b_{1}-s^{2}V}{V}\Phi_{x}+\frac{R}{V}\widehat{W}_{x}
-\frac{\mu}{V}\Psi_{xx}+\frac{\gamma-1}{V}U_{x}\Psi=F_{1},}\\
{\displaystyle\frac{R}{\gamma-1}\widehat{W}_{t}+(b_{1}-s^{2}V)\Psi_{x}-
\frac{\kappa}{V}\left(\widehat{W}_{x}+\frac{\gamma-1}{R}U_{x}\Psi\right)_{x}}\\
{\displaystyle\quad-s
U_{x}\Psi+\frac{\kappa}{V^{2}}\Theta_{x}\Phi_{x}=F_{2},}
\end{array}\right.
\eqno(3.8)
$$
where $F_{1}$ and $F_{2}$ are nonlinear terms with respect to
$(\Phi,\Psi,\widehat{W})$,  that is
$$
\left\{
\begin{array}{lll}
{\displaystyle
F_{1}=\frac{\gamma-1}{2V}\psi^{2}-\frac{\phi}{V(V+\phi)}\left\{(b_{1}-s^{2}V)\phi
-R w+\mu\psi_{x}\right\}},\\
{\displaystyle
F_{2}=-\frac{\kappa(\gamma-1)}{RV}\psi\psi_{x}+\frac{\psi}{V+\phi}
\left\{(b_{1}-s^{2}V)\phi-Rw+\mu
\psi_{x}\right\}}\\
{\displaystyle\qquad\quad
-\frac{\kappa\phi}{V(V+\phi)}\left(w_{x}-\frac{\Theta_{x}\phi}{V}\right).}
\end{array}\right.
\eqno(3.9)
$$

By (3.3)-(3.4), the initial values satisfy
$$
\begin{array}{lll}
\Phi(x,0)&{\displaystyle=-\int_{x}^{+\infty}[v_{0}(y)-V(y+\alpha-\beta)]dy}\\
&={\displaystyle\tilde{\Phi}_{0}(x)+\int_{x}^{+\infty}[V(y+\alpha-\beta)-V(y-\beta)]dy}\\
&={\displaystyle\tilde{\Phi}_{0}(x)+\int_{0}^{\alpha}[v_{+}-V(x+\varsigma-\beta)]
d\varsigma=:\Phi_{0}(x).}
\end{array}
\eqno(3.10)
$$
$$
\begin{array}{lll}
\Psi(x,0)&{\displaystyle=-\int_{x}^{+\infty}[u_{0}(y)-U(y+\alpha-\beta)]dy}\\
&={\displaystyle\widetilde{\Psi}_{0}(x)+\int_{0}^{\alpha}[u_{+}-U(x+\varsigma-\beta)]
d\varsigma=:\Psi_{0}(x)}.
\end{array}
\eqno(3.11)
$$
$$
\begin{array}{lll}
W(x,0){\displaystyle=-\int_{x}^{+\infty}
\left[(\frac{R\theta_{0}}{\gamma-1}+\frac{u_{0}^{2}}{2})(y)-
(\frac{R\Theta}{\gamma-1}+\frac{U^{2}}{2})(y+\alpha-\beta)\right]dy}\\[3mm]
={\displaystyle\widetilde{W}_{0}(x)+\int_{x}^{+\infty}
\left[(\frac{R\Theta}{\gamma-1}+\frac{U^{2}}{2})(y+\alpha-\beta)
-(\frac{R\Theta}{\gamma-1}+\frac{U^{2}}{2})(y)\right]dy}\\[3mm]
={\displaystyle\widetilde{W}_{0}(x)+
\int_{0}^{\alpha}\frac{R}{\gamma-1}[\theta_{+}-\Theta(x+\varsigma
-\beta)]+\frac{1}{2}[u_{+}^{2}-U^{2}(x+\varsigma -\beta)]
d\varsigma}\\=:W_{0}(x).
\end{array}
\eqno(3.12)
$$
$$
\widehat{W}(x,0)=\frac{\gamma-1}{R}[W_0(x)-U(x+\alpha-\beta)\Psi_0(x)]=:\widehat{W}_{0}(x).
\eqno(3.13)
$$
Furthermore, by the same way as in \cite{[9]}, we have

\noindent \textbf{Lemma 3.1.}\textit{ Under the assumptions}
(2.11)\textit{and} (2.13),
$(\widetilde{\Phi}_{0},\widetilde{\Psi}_{0},\widetilde{W}_{0})\in
H^{2} $ \textit{and the shift}
$$\alpha\rightarrow0
\quad \mathrm{as}\
 \|(\widetilde{\Phi}_{0},\widetilde{\Psi}_{0},\widetilde{W}_{0})\|_{2}\rightarrow0\
 \mathrm{and }\ \beta\rightarrow+\infty.\eqno(3.14)$$

\

\noindent \textbf{Lemma 3.2.} \textit{ Under the assumptions} (2.11)
\textit{and} (2.13), \textit{the initial perturbations}
$(\Phi_{0},\Psi_{0},\widehat{W}_{0})\in H^{2}$\textit{and satisfies}
$$\|(\Phi_{0},\Psi_{0},\widehat{W}_{0})\|\rightarrow0\quad \mathrm{as}\
\|(\widetilde{\Phi}_{0},\widetilde{\Psi}_{0},\widetilde{W}_{0})\|\rightarrow
0\ \mathrm{and} \ \beta\rightarrow +\infty.\eqno (3.15)$$

\

By (3.3) (3.5) and (2.5), the boundary values satisfy
$$
\begin{array}{lll}
\displaystyle
\Psi(0,t)&=&\displaystyle -\int_{0}^{+\infty}\psi(y,t)dy\\
\displaystyle&=&\displaystyle -s\int_{t}^{+\infty}\left[ U(-s
\tau+\alpha-\beta)-u_{-}\right]d\tau:=A(t),
\end{array}
\eqno(3.16)$$
$$
\widehat{W}_{x}(0,t)-\frac{\gamma-1}{2R}\Psi_x^2(0,t)
=\omega(0,t)-U_x(-st+\alpha-\beta)A(t):=B(t). \eqno(3.17)
$$

For any $T>0$, we define the solution space of the problem (3.5),
with the initial values (3.10), (3.11), (3.13) and the boundary
values (3.16), (3.17) by
$$
X_{m,M}(0,T)= \left\{
\begin{array}{l}
\displaystyle (\Phi,\Psi,\widehat{W}):\ (\Phi,\Psi,\widehat{W})\in
C(0,T;H^{2});\\
\displaystyle\  \Phi_{x}\in
 L^{2}(0,T;H^{1});\ (\Psi_{x};\widehat{W}_{x})\in
L^{2}(0,T;H^{2});\\
\displaystyle\ \sup_{t\in[0,T]}\|(\Phi,\Psi,W)(t)\|_{2}\leq M; \
\inf_{x,t}(V+\Phi_x)\geq m
\end{array}
\right\}
 \eqno(3.18)
$$
where $T,M,m$ are the positive constants.

\section{Proof of Theorem 2.1}\label{sec:int}

In this section, we give the proof of the Theorem 2.1. Without
loss of generality, we may restrict $\beta>1$ and $|\alpha|<1$.
First we state the local existence result for the IBVP (3.8),
(3.10)-(3.13) and (3.16)-(3.17), whose proof is given in section
5.

\

\noindent\textbf{Proposition 4.1.}(Local Existence) \textit{There
exists a positive constant $b$ such that if
$\|(\Phi_0,\Psi_0,\widehat{W}_0)\|_2\leq M$, and if
$\inf_{x,t}(V+\Phi_{0x})\geq m>0$,then there exists a positive
constant $T_0=T_0(m,M)$ such that the system (3.8), with the
initial values (3.10), (3.11), (3.13) and the boundary values
(3.16), (3.17), has a unique solution $(\Phi,\Psi,\widehat{W})\in
X_{\frac12m,bM}(0,T_0)$.}

\

Denote that
$$
\begin{array}{lll} {\displaystyle
N(T)=\sup_{\tau\in[0,T]}(\|\Phi(\tau)\|_{2}+\|\Psi(\tau)\|_{2}+\|W(\tau)\|_{2}),}\\
N_{0}=\|\Phi_{0}\|_{2}+\|\Psi_{0}\|_{2}+\|W_{0}\|_{2}.
\end{array}
$$

\

\noindent\textbf{Proposition 4.2.}(A Priori Estimates)
\textit{Let} $(\Phi,\Psi,W)\in X_{\frac12m,b\varepsilon}(0,T)$
\textit{be a solution of the problem (3.5) and} $1<\gamma\leq2$.
\textit{Then there exist positive constants}
$\delta_1,\varepsilon_1$ \textit{and} $C$,\textit{ which are
independent of} $T$, \textit{such that if}
$(\gamma-1)d\leq\delta_1$ \textit{and}
$N_0+\varepsilon+\beta^{-1}\leq\varepsilon_1$,\textit{ then the
following estimate holds for} $t\in [0,T]$
$$
\begin{array}{lll}
{\displaystyle
\|(\Phi,\Psi,\frac{\widehat{W}}{(\gamma-1)^{\frac{1}{2}}})(t)\|^{2}+\|(\phi,\psi,\frac{w}{(\gamma-1)^{\frac{1}{2}}})(t)\|_{1}^{2}
+\int_{0}^{t}\|(\psi,w)\|^{2}_{2}+\|\phi\|_{1}^{2}d\tau}\\
{\displaystyle\leq C\left (N_0+e^{-cd\beta}\right)}.
\end{array}
\eqno(4.1)
$$

 \

With the local existence Proposition 4.1 in hand, for the proof of
the Theorem 2.1 by the standard continuum process, it is sufficient
to prove the a priori estimate Proposition 4.2. In order to prove
the Proposition 4.2, we first give some Lemmas. The following Lemma
is about the boundary estimates.

\

\noindent\textbf{Lemma 4.3.} \textit{For} $0\leq t\leq T$,
\textit{the following inequalities hold}:
$$
\begin{array}{l}
\displaystyle\int_{0}^{t}(\Phi\Psi)\bigl|_{x=0}d\tau,
\int_{0}^{t}(\Psi\Psi_{x})\bigl|_{x=0}d\tau,\
\int_{0}^{t}(\widehat{W}\Psi) \bigl|_{x=0}d\tau, \int_{0}^{t}(\psi
w)\bigl|_{x=0}d\tau \leq
Ce^{-cd\beta},\\[3mm]
\displaystyle
\int_{0}^{t}(\widehat{W}_{x}\widehat{W})\bigl|_{x=0}d\tau \leq
Ce^{-cd\beta}+CN(T)\int_{0}^{t}(\|\Psi_{x}\|^{2}+\|\Psi_{xx}\|^2)d\tau,\\[3mm]
\displaystyle \int_{0}^{t}(\phi\psi)\bigl|_{x=0}d\tau,\
 \int_{0}^{t}(\psi\psi_{x})\bigl|_{x=0}d\tau,\
 \int_{0}^{t}(\psi_x\psi_\tau)\bigl|_{x=0}d\tau \leq
C(e^{-cd\beta}+\|\phi_0\|_{1}),\\
\displaystyle\int_0^t(ww_x)\big|_{x=0}d\tau,\
\int_0^t(w_xw_\tau)\big|_{x=0}d\tau\leq
(\gamma-1)d\int_0^t\|w_{xx}\|^2(\tau)d\tau+ Ce^{-cd\beta}.
\end{array}
$$

\

 \noindent\textbf{Proof.}
Since $s>0$, and $\beta\gg1,
 |\alpha|<1$, we have from (2.3) and (3.16) that
 $$
 |\Psi(0,t)|=|A(t)|\leq Ce^{-cd\beta}e^{-cd
 t}.
 $$
Thus,
$$
\int_{0}^{t}(\Phi\Psi)\bigl|_{x=0}d\tau\leq
Cd^{-1}N(T)e^{-cd\beta}\leq Ce^{-cd\beta}.
$$
Similarly we can estimate the term
$\displaystyle\int_{0}^{t}(\Psi\Psi_{x})\bigl|_{x=0}d\tau,\
\int_{0}^{t}(\widehat{W}\Psi) \bigl|_{x=0}d\tau$.

Also,
$$
\int_0^t(\psi w)\big|_{x=0}d\tau\leq
N(T)\int_0^t|\theta_{-}-\Theta(-s\tau+\alpha-\beta)|d\tau\leq
Ce^{-cd\beta}.
$$
From (3.8),
$$\widehat{W}_{x}(0,t)=w(0,t)-\frac{\gamma-1}{R}(U_{x}\Psi(0,t)-\frac{\Psi_{x}^{2}(0,t)}{2}),$$
so we have from (2.3) that
$$
\int_{0}^{t}(\widehat{W}_{x}\widehat{W})\bigl|_{x=0}d\tau \leq
Ce^{-cd\beta}+CN(T)\int_{0}^{t}(\|\Psi_{x}\|^{2}+\|\Psi_{xx}\|^2)d\tau.
$$
By using the free boundary condition in (1.5), one has
$$\frac{R\theta_{-}}{v(0,t)}-\mu\frac{v(0,t)_{t}}{v(0,t)}=\frac{R\theta_{-}}{v_{-}},$$
and then
$$
\begin{array}{ll}
\di v(0,t)-v_{-}&\di =(v_{0}(0)-v_{-})e^{-\frac{p_{0}}{\mu}t}\\
\di&=(V(\a-\b)-v_-+\phi_0(0))e^{-\frac{p_{0}}{\mu}t}\\
\di &\leq C(e^{-cd\beta}+\|\phi_0\|_{1})e^{-\frac{p_{0}}{\mu}t}
\end{array}
 \eqno(4.2)
$$

By using (2.3) and (4.2), we obtain
$$
\begin{array}{ll}
\di |\phi(0,t)|&\di =|v(0,t)-V(-st+\alpha-\beta)|\\
\di&\leq |v(0,t)-v_-|+|V(-st+\alpha-\beta)-v_-| \\
 \di &\leq
C(e^{-cd\beta}+\|\phi_0\|_{1})e^{-\frac{p_{0}}{\mu}t}+
Ce^{-cd\beta}e^{-cd t},
\end{array}
$$
$$
\begin{array}{ll}
\di |\psi_{x}(0,t)|&\di =
 |\frac{p_0}{\mu}(v_{-}-v_{0}(0))e^{-\frac{p_{0}}{\mu}t}-U_{x}(-st+\alpha-\beta)|\\
\di& \leq C(e^{-cd\beta}+\|\phi_0\|_{1})e^{-\frac{p_{0}}{\mu}t}+
Ce^{-cd\beta}e^{-cd t}.
\end{array}
$$
Then we get at once
$$
\begin{array}{ll}
\di\int_{0}^{t}(\phi\psi)\bigl|_{x=0}d\tau,\quad\int_{0}^{t}(\psi\psi_{x})\bigl|_{x=0}d\tau
&\di\leq CN(T)(e^{-cd\beta}+\|\phi_0\|_{1})\\
&\di\leq C(e^{-cd\beta}+\|\phi_0\|_{1}),
\end{array}
$$ and
$$
\begin{array}{lll}
\displaystyle\int_{0}^{t}(\psi_x\psi_\tau)\bigl|_{x=0}d\tau
&=&\displaystyle\int_{0}^{t}\left(\psi_x(0,\tau)\psi(0,\tau)\right)_\tau
d\tau
-\int_{0}^{t}(\psi_{x\tau}\psi)\bigl|_{x=0}d\tau\\
&=&\displaystyle
\psi_x(0,\tau)\psi(0,\tau)\big|_0^t-\int_{0}^{t}(\psi_{x\tau}\psi)\bigl|_{x=0}d\tau\\
&\leq&\displaystyle CN(T)(e^{-cd\beta}+\|\phi_0\|_{1})\leq C
(e^{-cd\beta}+\|\phi_0\|_{1}).
\end{array}
$$
Finally,
$$
\begin{array}{ll}
&\displaystyle\int_0^t(ww_x)\big|_{x=0}d\tau \leq C(\gamma-1)d
e^{-cd\beta}\int_0^te^{-cd
\tau}\|w_x\|^{\frac12}\|w_{xx}\|^{\frac12}d\tau\\[3mm]
&\qquad\displaystyle \leq (\gamma-1)d
\int_0^t\|w_{xx}\|^2(\tau)d\tau+C(\gamma-1)d e^{-cd\beta}\int_0^t
\|w_x\|^{\frac23}(\tau)e^{-{\frac23}cd \tau}d\tau\\[3mm]
&\qquad\displaystyle \leq (\gamma-1)d
\int_0^t\|w_{xx}\|^2(\tau)d\tau+C(\gamma-1)d
 e^{-cd\beta}[N(T)]^{\frac23}\int_0^t
e^{-{\frac23}cd \tau}d\tau\\[3mm]
&\qquad\displaystyle \leq (\gamma-1)d
\int_0^t\|w_{xx}\|^2(\tau)d\tau+Ce^{-cd\beta},
\end{array}
$$
and
$$
\begin{array}{ll}
&\displaystyle \int_0^t(w_xw_\tau)\big|_{x=0}d\tau\leq
C(\gamma-1)d^2e^{-cd\beta}\int_0^te^{-cd
\tau}\|w_x\|^{\frac12}\|w_{xx}\|^{\frac12}d\tau\\[3mm]
&\qquad\displaystyle \leq
(\gamma-1)d\int_0^t\|w_{xx}\|^2(\tau)d\tau+C(\gamma-1)d^2e^{-cd\beta}\int_0^t
\|w_x\|^{\frac23}(\tau)e^{-{\frac23}cd \tau}d\tau\\[3mm]
&\qquad\displaystyle \leq
(\gamma-1)d\int_0^t\|w_{xx}\|^2(\tau)d\tau+C(\gamma-1)d^2e^{-cd\beta}[N(T)]^{\frac23}\int_0^t
e^{-{\frac23}cd \tau}d\tau\\[3mm]
&\qquad\displaystyle \leq
(\gamma-1)d\int_0^t\|w_{xx}\|^2(\tau)d\tau+Ce^{-cd\beta}.
\end{array}
$$
 We complete the proof of the lemma 4.3.

\

\textbf{Lemma 4.4.} \textit{For} $(\gamma-1)d\leq\delta_0$
\textit{small enough, then}
$$
\begin{array}{lll}
{\displaystyle
\|(\Phi,\Psi,\frac{\widehat{W}}{(\gamma-1)^{\frac{1}{2}}})(t)\|^{2}+
\int_{0}^{t}\|\;|V_{x}|^{\frac{1}{2}}(\Psi,\frac{\widehat{W}}{(\gamma-1)^{\frac{1}{2}}})
(\tau)\|^{2}d\tau}\\
{\displaystyle+
\int_{0}^{t}\|(\Psi_{x},\widehat{W}_{x})(\tau)\|^{2}d\tau-C(\gamma-1)d\int_{0}^{t}\|\Phi_{x}(\tau)\|^{2}d\tau}\\
{\displaystyle\leq
C\left\{\|(\Phi_{0},\Psi_{0},\frac{\widehat{W}_{0}}{(\gamma-1)^{\frac{1}{2}}})\|^{2}
+\int_{0}^{t}\int_0^{+\infty}|\Psi|\:|F_{1}|+|\widehat{W}|\:|F_{2}| dxd\tau\right\}}\\
{\displaystyle\quad+ CN(T)\int_{0}^{t}\|\Psi_{xx}\|^{2}d\tau
+Ce^{-cd\beta}.}
\end{array}
\eqno(4.3)
$$

\

\textbf{Proof.} Let
$$
k(V)=(b_{1}-s^{2}V)^{-1}.
$$ Multiplying the
first equation of (3.8) by $\Phi$, the second equation of (3.8) by
$k(V)V\Psi$ and the third equation of (3.8) by
$Rk(V)^{2}\widehat{W}$, respectively, summing them up, we have
$$
\begin{array}{lll}
{\displaystyle
E_{1}(\Phi,\Psi,\widehat{W})_{t}+E_{2}(\Psi,\Psi_{x})+
E_{3}(\widehat{W},\widehat{W}_{x})+G(\Psi,\widehat{W},\Phi_{x},\widehat{W}_{x})}\\
{\displaystyle+\left\{\mu k(V)\Psi \Psi_{x}-\Phi\Psi-\frac{R\kappa
k(V)^{2}}{V}(\widehat{W}_{x}+\frac{\gamma-1}{R}U_{x}\Psi)\widehat{W}+Rk(V)\widehat{W}\Psi\right\}_{x}}\\
{\displaystyle=k(V)V\Psi F_{1}+Rk(V)^{2}\widehat{W}F_{2},}
\end{array}
\eqno(4.4)
$$
where
$$
\begin{array}{lll}
{\displaystyle E_{1}(\Phi,\Psi,\widehat{W})=\frac{1}{2}\left(\Phi^{2}+
k(V)V\Psi^{2}+\frac{R^{2}}{\gamma-1}k(V)^{2}\widehat{W}^{2}\right),}\\
{\displaystyle
E_{2}(\Psi,\Psi_{x})=\left[\frac{s}{2}(k(V)V)_{x}+(\gamma-1)k(V)U_{x}\right]\Psi^{2}
+\mu k(V)\Psi_{x}^{2}+\mu k(V)_{x}\Psi\Psi_{x},}\\
{\displaystyle E_{3}(\widehat{W},\widehat{W}_{x})=\frac{s
R^{2}}{\gamma-1}k(V)k(V)_{x}\widehat{W}^{2}+\kappa
R\frac{k(V)^{2}}{V}\widehat{W}_{x}^{2}),} \\
{\displaystyle G(\Psi,\widehat{W},\Phi_{x},\widehat{W}_{x})=\kappa
R\left(\frac{k(V)^{2}}{V}\right)_{x}\widehat{W}
\left(\widehat{W}_{x}+\frac{\gamma-1}{R}U_{x}\Psi\right)}\\
{\displaystyle\qquad\qquad\qquad\qquad+\kappa
R\frac{k(V)^{2}}{V^{2}}\Phi_{x}\Theta_{x}\widehat{W}+
\kappa(\gamma-1)\frac{k(V)^{2}}{V}U_{x}\Psi
\widehat{W}_{x},}
\end{array}
$$
Since
$$
p_{-}\leq k(V)^{-1}=b_{1}-s^{2}V\leq p_{+}, \eqno(4.5)
$$
one has
$$
c\left(\Phi^{2}+\Psi^{2}+\frac{\widehat{W}^{2}}{\gamma-1}\right)\leq
E_{1}\leq C
\left(\Phi^{2}+\Psi^{2}+\frac{\widehat{W}^{2}}{\gamma-1}\right),
\eqno(4.6)
$$
$$
E_{3}\geq
c\left(|V_{x}|\frac{\widehat{W}^{2}}{\gamma-1}+\widehat{W}_{x}^{2}\right),
\eqno(4.7)
$$
and for $\forall \alpha_1>0$, there $\exists$ a constant
$C_{\alpha_1}$ such that
$$
|G|\leq\alpha_1\left(|V_{x}|\frac{\widehat{W}^{2}}{\gamma-1}+\widehat{W}_{x}^{2}\right)+C_{\alpha_1}(\gamma-1)
d\left[|V_{x}|\left(\Psi^{2}+\frac{\widehat{W}^{2}}{\gamma-1}\right)+\Phi_{x}^{2}\right].
\eqno(4.8)
$$
By using the method in \cite{[7]}, for $\gamma\in (1,2]$ and
suitably small $(\gamma-1)d>0$, one has
$$
\begin{array}{lll}
{\displaystyle
\inf_{x>0}\frac{\frac{s}{2}(k(V)V)_{x}+(\gamma-1)k(V)U_{x}}{V_{x}}>0,}\\
{\displaystyle \sup_{x>0} \frac{\mu\left\{\mu|k(V)_{x}|^{2}-
4\left[\frac{s}{2}(k(V)V)_{x}+(\gamma-1)k(V)U_{x}\right]k(V)\right\}}{V_{x}}<0,}\end{array}
\eqno(4.9)
$$
and then we get
$$
E_{2}\geq c(|V_{x}\Psi^{2}|+\Psi_{x}^{2}).
\eqno(4.10)
$$
Combining with the boundary estimates in Lemma 4.3, (4.3) is
obtained.

\

\textbf{Lemma 4.5.} \textit{There is a constant} $C$ \textit{such
that}
$$
\begin{array}{lll}
{\displaystyle\|\phi(t)\|^{2}+\int_{0}^{t}\|\phi(\tau)\|^{2}d\tau}\\
{\displaystyle-C\{\|\Psi(t)\|^{2}+
\int_{0}^{t}\|\;|V_{x}|^{\frac{1}{2}}\Psi(\tau)\|^{2}+\|(\psi,\widehat{W}_{x})(\tau)\|^{2}d\tau\}}\\
{\displaystyle\leq
C\left\{\|\Psi_{0}\|^{2}+\|\phi_{0}\|^{2}+\int_{0}^{t}\int_0^{+\infty}|\phi|\;|F_{1}|dxd\tau\right\}.}
\end{array}
\eqno(4.11)
$$

\

\textbf{Proof.} Multiplying $(3.8)_{1}$ by $V\Psi_x-V_x\Psi$,
$(3.8)_{2}$ by $-V\Phi_x$, then applying $\partial_x$ to $(3.8)_{1}$
and multiplying the resulting equation by $\mu\Phi_x$, calculating
all their sums, we get

$$
\begin{array}{l}
\displaystyle(\frac{\mu\Phi_x^2}{2}-V\Phi_x\Psi)_t+(V\Psi\Psi_x)_x+(b_1-s^2V)\Phi_x^2\\
\displaystyle\qquad=V_x\Psi\Psi_x+V\Psi_x^2+\left[R\widehat{W}_x-s(\gamma-1)V_x\Psi-V_t\Psi-VF_1\right]\Phi_{x}.
\end{array}
\eqno(4.12)
$$
Integrating (4.12) over $[0,+\infty)\times[0,t]$ with respect to
$x,t$ and using the boundary estimates Lemma 4.2, we obtain Lemma
4.5.

\

From (3.7), we easily have
$$\begin{array}{lll}
{\displaystyle \int_{0}^{t}\|w(\tau)\|^{2}d\tau-
C\left\{\int_{0}^{t}\|\;|V_{x}|^{\frac{1}{2}}\Psi(\tau)\|^{2}+\|\widehat{W}_{x}(\tau)\|^{2}d\tau\right\}}\\
{\displaystyle\qquad\leq
C\int_{0}^{t}|\psi^{2}w|dxd\tau.}\end{array} \eqno(4.13)$$

\

Now we rewrite (3.2) in the form

$$
\left\{
\begin{array}{llll}
{\displaystyle \phi_{t}-\psi_{x}=0,}\\
{\displaystyle
\psi_{t}-\frac{b_{1}-s^{2}V}{V}\phi_{x}+\frac{R}{V}w_{x}-(\frac{\mu}{V}\psi_{x})_{x}}
{\displaystyle-\{\frac{b_{1}-s^{2}V}{V}\}_{x}\phi+(\frac{R}{V})_{x}w}=f_{1}\\
{\displaystyle\frac{R}{\gamma-1}w_{t}+(b_{1}-s^{2}V)\psi_{x}-(\frac{\kappa}{V}w_{x})_{x}
+(\frac{\kappa}{V^{2}}\Theta_{x}\phi)_{x}}\\
{\displaystyle\quad-\frac{1}{V}\{(b_{1}-s^{2}V)\phi-Rw+\mu\psi_{x}\}U_{x}=f_{2},}
\end{array}
\right.
\eqno(4.14)
$$
where $f_{1}$ and $f_{2}$ are nonlinear terms with respect to
$(\phi,\psi,w)$

$$\begin{array}{llll}
{\displaystyle
f_{1}=-\{\frac{\phi}{V(V+\phi)}[(b_{1}-s^{2}V)\phi-Rw+\mu
\psi_{x}]\}_{x},}\\
{\displaystyle
f_{2}=\frac{1}{V+\phi}\{(b_{1}-s^{2}V)\phi-Rw+\mu\psi_{x}\}(\psi_{x}-\frac{1}{V}U_{x}\phi)}\\
{\displaystyle\qquad-\{\frac{\kappa \phi
}{V(V+\phi)}(w_{x}-\frac{1}{V}\Theta_{x}\phi)\}_{x}.}\end{array}
$$
The following Lemma is the estimates of $(\phi,\psi,\omega)$ and
$(\phi_x,\psi_x,\omega_x)$.

\

 \textbf{ Lemma 4.6.} \textit{There is a constant} $C$ \textit{such that}

 $$\begin{array}{llll}
{\displaystyle
\|(\phi,\psi,\frac{w}{(\gamma-1)^{\frac{1}{2}}})(t)\|^{2}+\int_{0}^{t}\|\partial_{x}(\psi,w)(\tau)
 \|^2d\tau-C\int_{0}^{t}\|(\phi,\psi,w)(\tau)\|^{2}d\tau}\\
 {\displaystyle\qquad-(\gamma-1)d\int_0^t\|w_{xx}\|^2(\tau)d\tau\leq C\|(\phi_{0},\psi_{0},\frac{w_{0}}
 {(\gamma-1)^{\frac{1}{2}}})\|^{2}}\\
 {\displaystyle\qquad+C\int_{0}^{t}\int_{0}^{+\infty}(|\psi||f_{1}|+|w||f_{2}|)dxd\tau+C(e^{-cd\beta}+\|\phi_0\|_1).}\\
 {\displaystyle\|\phi_{x}(t)\|^{2}+\int_{0}^{t}\|\phi_{x}(\tau)\|^{2}d\tau-C\left\{\|\psi(t)\|^{2}+
 \int_{0}^{t}\|(\psi,w)(\tau)\|_{1}^{2}d\tau\right\}}\\
  {\displaystyle\qquad\leq C\left\{\|\psi_{0}\|^{2}+
  \int_{0}^{t}\int_{0}^{+\infty}|\phi_{x}||f_{1}|dxd\tau+(e^{-cd\beta}+\|\phi_{0}\|_1)\right\}.}\\
 {\displaystyle
\|\partial_{x}(\psi,\frac{w}{(\gamma-1)^{\frac{1}{2}}})(t)\|^{2}+\int_{0}^{t}\|\partial_{xx}(\psi,w)(\tau)
 \|^2d\tau-C\int_{0}^{t}\|(\phi,\psi,w)(\tau)\|_1^{2}d\tau}\\
 {\displaystyle\qquad\leq
 C\|\partial_{x}(\psi_{0},\frac{w_{0}}{(\gamma-1)^{\frac{1}{2}}})\|^{2}+C(e^{-cd\beta}}+\|\phi_0\|_1)\\
 {\displaystyle \qquad+ C\int_{0}^{t}\int_{0}^{+\infty}(|\psi_{xx}||f_{1}|+|w_{xx}||f_{2}|)dxd\tau.}
 \end{array}
\eqno(4.15)$$

\textbf{Proof.} Multiplying $(4.14)_1$ by $\phi$, $(4.14)_2$ by
$Vk(V)\psi$, $(4.14)_3$ by $Rk^2(V)w$,and adding them and
integrating over $x,t$, we have
$$
\begin{array}{l}
{\displaystyle\|(\phi,\psi,\frac{w}{(\gamma-1)^{\frac{1}{2}}})(t)\|^{2}+\int_{0}^{t}\|\partial_{x}(\psi,w)(\tau)
 \|^2d\tau-C\int_{0}^{t}\|(\phi,\psi,w)(\tau)\|^{2}d\tau}\\
 {\displaystyle+\int_0^t\left[-\phi\psi+Rk(V)\psi w-\mu k(V)\psi\psi_x+\frac{\kappa R}{V^2}\Theta_xk^2(V)\phi w
 -\frac{\kappa Rk^2(V)}{V}ww_x\right]_{x=0}d\tau}\\
 {\displaystyle\quad\leq C\|(\phi_{0},\psi_{0},\frac{w_{0}}
 {(\gamma-1)^{\frac{1}{2}}})\|^{2}+C\int_{0}^{t}\int_{0}^{+\infty}(|\psi||f_{1}|+|w||f_{2}|)dxd\tau.}
\end{array}
$$
Using the boundary estimate in Lemma 4.3, we can get the first
inequality in (4.15).

Now we want to get the estimate of $\|\phi_x\|^2$ in $(4.15)_2$.
Multiplying $(4.14)_{1}$ by $V\psi_x-V_x\psi$, $(4.14)_{2}$ by
$-V\phi_x$, then applying $\partial_x$ to $(4.14)_{1}$ and
multiplying the resulting equation by $\mu\phi_x$, calculating all
their sums, we get
$$
\begin{array}{l}
\displaystyle(\frac{\mu\phi_x^2}{2}-V\phi_x\psi)_t+(V\psi\psi_x)_x+(b_1-s^2V)\phi_x^2=V_x\psi\psi_x+V\psi_x^2+\\[3mm]
\displaystyle\quad\left[R\omega_x+V(\frac{\mu}{V})_x\psi_x+V(\frac{b_1-s^2V}{V})_x\phi-(\frac{R}{V})_x\omega-V_t\psi-Vf_1\right]\phi_{x}.
\end{array}
\eqno(4.16)
$$
Thus integrating the equation (4.16) and using the boundary estimate
Lemma 4.3,we obtain $(4.15)_2$.

 Multiplying $(4.14)_2$ by $-\psi_{xx}$, $(4.14)_3$ by
$-w_{xx}$ to get
$$
\begin{array}{ll}
\displaystyle(\frac12\psi_x^2+\frac{R}{2(\gamma-1)}w_x^2)_t-(\psi_x\psi_t+w_xw_t)_x+\frac{\mu}{V}\psi_{xx}^2+\frac{\kappa}{V}w_{xx}^2\\[3mm]
\displaystyle\quad=-\psi_{xx}\left[\frac{b_1-s^2V}{V}\phi_x-\frac{R}{V}w_x+(\frac{\mu}{V})_x\psi_x+(\frac{b_1-s^2V}{V})_x\phi-(\frac{R}{V})_xw+f_1\right]\\[3mm]
\displaystyle\quad=-w_{xx}\left[-(b_1-s^2V)\phi_x+(\frac{\kappa}{V})_xw_x+(\frac{\mu}{V^2}\Theta_x\phi)_x\right.\\[3mm]
\displaystyle\qquad\qquad\qquad\qquad\qquad\left.+\frac{1}{V}\{(b_1-s^2V)\phi-Rw+\mu\psi_x\}U_x+f_2\right]
\end{array}
$$
Integrating the above equality and using Lemma 4.3, we can get the
third inequality of (4.15). The proof of Lemma 4.6 is complete.

 \

Since
$$
\begin{array}{ll}
\displaystyle
|F_{1},F_{2}|&\displaystyle=O(1)\left[|(\phi,\psi,w)|^{2}+|(\phi,\psi)||(\psi_{x},w_{x})|\right],\\[2mm]
\displaystyle
|f_{1},f_{2}|&\displaystyle=O(1)\left[|(\phi,w)|^{2}+|(\phi,w)||(\phi_{x},\psi_{x},w_{x})|\right.\\[2mm]
&\displaystyle\quad\left.+|(\phi_{x},\psi_{x})||(\psi_{x},w_{x})|+|\phi||(\psi_{xx},w_{xx})|\right],
\end{array}
\eqno (4.17)
$$
combining (4.17) with the estimates (4.3), (4.11), (4.13), (4.15)
and using the a priori assumption $N(T)\leq b\v$ sufficiently small,
and also letting $(\gamma-1)d$ small enough, we can get the
following estimate
$$
N^2(T)+\int_{0}^{T}\|(\psi,w)\|^{2}_{2}+\|\phi\|_{1}^{2}d\tau\leq
\bar{C}(N_0+e^{-cd\beta}).
$$
where the constant $\bar{C}$ is independent of $T$. Thus we get the
desired a priori estimate (4.1) if we choose $N_0$ and
$e^{-cd\beta}$ small enough.

\section{The Local Existence}
In this section, we prove the local existence result Proposition
4.1 by the iteration method.  First we rewrite the equation (3.8)
with the initial values (3.10)-(3.13) and the boundary values
(3.16)-(3.17) as the following
$$
\left\{
\begin{array}{l}
\displaystyle\Psi_t-\frac{\mu}{V+\Phi_x}\Psi_{xx}=g_1:=g_1(\Psi,\Phi_x,\Psi_x,\widehat{W}_x),\\[3mm]
\displaystyle\Psi(0,t)=A(t),\\[2mm]
\displaystyle\Psi(x,0)=\Psi_0(x),
\end{array}
\right.  \eqno(5.1)
$$
$$
\left\{
\begin{array}{l}
\displaystyle
\frac{R}{\gamma-1}\widehat{W}_t-\frac{\kappa}{V+\Phi_x}(\widehat{W}_x-\frac{\gamma-1}{2R}\Psi_{x}^2)_x=g_2:=g_2(\Psi,\Phi_x,\Psi_x,\widehat{W}_x,\Psi_{xx}),\\[3mm]
\displaystyle\widehat{W}_x(0,t)-\frac{\gamma-1}{2R}\Psi_{x}^2(0,t)=B(t),\\[2mm]
\displaystyle\widehat{W}(x,0)=\widehat{W}_0(x),
\end{array}
\right.
\eqno(5.2)
$$
and
$$
\Phi(x,t)=\Phi_0(x)+\int_0^t\Psi_x(x,\tau)d\tau, \eqno(5.3)
$$
 where $A(t),B(t)$ is given in (3.16), (3.17) respectively, and
$$
g_1=\frac{b_1-s^2V}{V+\Phi_x}\Phi_x
-\frac{R}{V+\Phi_x}(\widehat{W}_x+\frac{\gamma-1}{R}U_x\Psi-\frac{\gamma-1}{2R}\Psi_x^2),
\eqno(5.4)
$$
$$
\begin{array}{ll}
g_2=&\displaystyle-\frac{b_1-s^2V}{V+\Phi_x}\Phi_x\Psi_x
+\frac{\kappa(\gamma-1)}{R(V+\Phi_x)}(U_x\Psi)_x+sU_x\Psi-\frac{\kappa\Theta_x\Phi_x}{V(V+\Phi_x)}\\
&\displaystyle+\frac{\mu\Psi_x\Psi_{xx}}{V+\Phi_x}-\frac{R\Psi_x}{V+\Phi_x}[\widehat{W}_x+\frac{\gamma-1}{R}(U_x\Psi-\frac{\Psi_x^2}{2})].
\end{array}
\eqno(5.5)
$$
To use the iteration method, we approximate the initial values
$(\Phi_0,\Psi_0,\widehat{W}_0)\in H^2(0,+\infty)$ by
$(\Phi_{0k},\Psi_{0k},\widehat{W}_{0k})\in H^5(0,+\infty)$ which
will be determined later. For fixed $k$, we define the sequence
$\{(\Phi_k^{(n)},\Psi_k^{(n)},\widehat{W}_k^{(n)})(x,t)\}_{n=1}^{\infty}$
by
$$
(\Phi_k^{(0)},\Psi_k^{(0)},\widehat{W}_k^{(0)})(x,t)=(\Phi_{0k},\Psi_{0k},\widehat{W}_{0k})(x),
\eqno(5.6)
$$
and if
$(\Phi_k^{(n-1)},\Psi_k^{(n-1)},\widehat{W}_k^{(n-1)})(x,t)$ is
given, then we define
$(\Phi_k^{(n)},\Psi_k^{(n)},\widehat{W}_k^{(n)})(x,t)$ as the
following
$$
\left\{
\begin{array}{l}
\displaystyle\Psi_{kt}^{(n)}-\frac{\mu}{V+\Phi_{kx}^{(n-1)}}\Psi_{kxx}^{(n)}=g_1^{(n-1)}:=g_1(\Psi_k^{(n-1)},\Phi_{kx}^{(n-1)},\Psi_{kx}^{(n-1)},\widehat{W}_{kx}^{(n-1)}),\\[3mm]
\displaystyle\Psi_{k}^{(n)}(0,t)=A(t),\\[2mm]
\displaystyle\Psi_{k}^{(n)}(x,0)=\Psi_{0k}(x),
\end{array}
\right.\eqno(5.7)
$$
$$
\left\{
\begin{array}{l}
\displaystyle
\frac{R}{\gamma-1}\widehat{W}_{kt}^{(n)}-\frac{\kappa}{V+\Phi_{kx}^{(n-1)}}(\widehat{W}_{kx}^{(n)}-\frac{\gamma-1}{2R}{\Psi_{kx}^{(n)}}^2)_x\\
\displaystyle\qquad=g_2^{(n-1)}:=g_2(\Psi_{k}^{(n)},\Phi_{kx}^{(n-1)},\Psi_{kx}^{(n)},\widehat{W}_{kx}^{(n-1)},\Psi_{kxx}^{(n)}),\\[3mm]
\displaystyle\widehat{W}_{kx}^{(n)}(0,t)-\frac{\gamma-1}{2R}{\Psi_{kx}^{(n)}}^2(0,t)=B(t),\\[2mm]
\displaystyle\widehat{W}_{k}^{(n)}(x,0)=\widehat{W}_{0k}(x),
\end{array}
\right. \qquad\qquad\eqno(5.8)
$$
and
$$
\Phi_{k}^{(n)}(x,t)=\Phi_{0k}(x)+\int_0^t\Psi_{kx}^{(n)}(x,\tau)d\tau.
\eqno(5.9)
$$
Now we construct the approximate initial values
$(\Phi_{0k},\Psi_{0k},\widehat{W}_{0k})(x)$. Firstly we choose
$\Phi_{0k}\in H^5$ such that $\Phi_{0k}\rightarrow\Phi_0$ strongly
in $H^2$ as  $k\rightarrow\infty.$ Let
$$
\overline\Psi_{0}(x):=\Psi_0(x)-A(0)e^{-x^2}.
$$
Note that $A(0)=\Psi_0(0).$ Then we have $\overline\Psi_{0}(x)\in
H_0^2$. Now we choose $\overline\Psi_{0k}(x)\in H_0^3\cap H^5$
such that $\overline\Psi_{0k}\rightarrow\overline\Psi_{0}$
strongly in $H^2$ as $k\rightarrow\infty$. We construct
$$
\Psi_{0k}(x):=\overline\Psi_{0k}(x)+A(0)e^{-x^2}, \eqno(5.10)
$$
then we have
$\Psi_{0k}\rightarrow\overline\Psi_0(x)+A(0)e^{-x^2}=\Psi_0(x)$
strongly in $H^2$ as $k\rightarrow\infty$. Moreover,
$\Psi_{0k}(x)$ constructed in (5.10) satisfies the compatibility
condition $\Psi_{0k}(0)=A(0)$ for the approximate equation (5.7).
Now we turn to the compatibility condition for the equation (5.8).
Let
$$
\overline{\widehat{W}}_0(x):=\widehat{W}_0(x)-B(0)xe^{-x^2}-\widehat{W}_0(0)e^{-x^2}.
$$
It is obvious that $\overline{\widehat{W}}_0(x)\in H^2_0$. So we
can choose $\overline{\widehat{W}}_{0k}(x)\in H_0^3\cap H^5$ such
that
$\overline{\widehat{W}}_{0k}(x)\rightarrow\overline{\widehat{W}}_{0}(x)$
strongly in $H^2$ as $k\rightarrow\infty$. Set
$$
\widehat{W}_{0k}(x):=\overline{\widehat{W}}_{0k}(x)+B(0)xe^{-x^2}+\widehat{W}_{0}(0)e^{-x^2}.
\eqno(5.11)
$$
Then we have
$\widehat{W}_{0k}(x)\rightarrow\overline{\widehat{W}}_{0}(x)+B(0)xe^{-x^2}+\widehat{W}_{0}(0)e^{-x^2}=\widehat{W}_{0}(x)$
strongly in $H^2$ as $k\rightarrow\infty$. Note that
$B(0)=\widehat{W}_{0x}(0)-\frac{\gamma-1}{2R}\Psi^2_{0x}(0)$. We
verify that the approximated initial values
$\Psi_{0k}(x),\widehat{W}_{0k}(x)$ satisfy  the following
compatibility condition for the equation (5.8),
$$
\widehat{W}_{0kx}(0)-\frac{\gamma-1}{2R}\Psi_{0kx}^2(0)=\overline{\widehat{W}}_{0kx}(0)+B(0)-\frac{\gamma-1}{2R}\overline\Psi_{0kx}^2(0)=B(0).
$$
And it is easy to choose that the above approximation
$(\Phi_{0k}(x),\Psi_{0k}(x),\widehat{W}_{0k}(x))$ satisfies
$\|(\Phi_{0k},\Psi_{0k},\widehat{W}_{0k})\|_2\leq \frac32 M$ and
$\inf_{x}(V+\Phi_{0kx})\geq \frac23m$ for any fixed $k$.

If $(\Psi_k^{(n-1)},\Psi_k^{(n-1)},\widehat{W}_k^{(n-1)})\in
X_{\frac12m,bM}(0,t_0)\cap C(0,t_0;H^5)$, then $g_1^{(n-1)}\in
C(0,t_0;H^4).$ By linear parabolic theory, since $\Psi_{0k}\in
H^5$, there exists a unique solution to (5.7) satisfying
$$
\Psi_k^{(n)}\in C(0,t_0;H^5)\cap C^1(0,t_0;H^3)\cap
L^2(0,t_0;H^6).
$$
Substituting $\Psi_k^{(n)}$ into $g_2^{(n-1)}$, we have that
$g_2^{(n-1)}\in C(0,t_0;H^3)$. Using linear parabolic theory
again, we obtain
$$
\widehat{W}_k^{(n)}\in C(0,t_0;H^5)\cap C^1(0,t_0;H^3)\cap
L^2(0,T;H^6).
$$
 From (5.9), we also have
$$
\Phi_k^{(n)}\in C(0,t_0;H^5)\cap C^1(0,t_0;H^3)\cap
L^2(0,t_0;H^6).
$$
The elementary energy estimates to the equation (5.7)-(5.8) yield
that
$$
\|(\Psi_k^{(n)},\widehat{W}_k^{(n)})\|^2_2\leq (bM)^2.
$$
if the time interval $t_0=t_0(m,M)$ is suitably small. We omit the
detailed calculations for brevity.

Now from (5.9), we can compute that
$$
\|\Phi_k^{(n)}\|^2_2\leq (bM)^2,
$$
and
$$
\inf_{x,t\in[0,t_0]}(V+\Phi_{kx}^{(n)})\geq \frac12 m.
$$
Therefore we have
$(\Psi_k^{(n)},\Psi_k^{(n)},\widehat{W}_k^{(n)})\in
X_{\frac12m,bM}(0,t_0)\cap C(0,t_0;H^5)$. Since
$\|(\Psi_k^{(0)},\Psi_k^{(0)},\widehat{W}_k^{(0)})\|_5$ is
uniformly bounded for fixed $k$, we can show that $(\Psi_k^{(n)},$
$\Psi_k^{(n)},\widehat{W}_k^{(n)})$ is the Cauchy sequence in
$C(0,t_0;H^4)$. Letting $n\rightarrow\infty$ in (5.7)-(5.9), we
get a solution $(\Phi_{k},\Psi_{k},\widehat{W}_{k})(x,t)$ of
(5.1)-(5.3) with the initial values replaced by
$(\Phi_{0k},\Psi_{0k},\widehat{W}_{0k})(x)$ in the time interval
$[0,t_0]$.

In the same way we can show that
$(\Phi_{k},\Psi_{k},\widehat{W}_{k})(x)$ is a Cauchy sequence in
$C(0,T_0;H^2)$ (takin $T_0$ smaller than $t_0$ if necessary). Now
letting $k\rightarrow\infty$, we get the desired unique solution
$(\Phi,\Psi,\widehat{W})(x,t)$ to (5.1)-(5.3) in the time interval
$[0,T_0]$.

\vspace{1cm}\noindent {\bf Acknowledgements:}\,\, The research of F.
M. Huang was supported in part by NSFC Grant No. 10825102 for
distinguished youth scholar, NSFC-NSAF Grant No. 10676037 and 973
project of China, Grant No.2006CB805902. The research of Y. Wang was
supported by the NSFC grant (No. 10801128).

\end{document}